\documentclass{article}
\usepackage{amsfonts}
\usepackage{amsmath}
\usepackage[onehalfspacing]{setspace}
\usepackage[hyphens,spaces,T1]{url}

\newtheorem{theorem}{Theorem}

\newtheorem{lemma}[theorem]{Lemma}

\input{tcilatex}
\begin{document}

\begin{center}
{\LARGE Modeling Spatial Overdispersion with the Generalized Waring Process}

\bigskip

\textbf{Mimoza Zografi and Evdokia Xekalaki}

Department of Statistics, Athens University of Economics, Greece
\end{center}

\bigskip \textbf{Abstract}

\bigskip Modeling spatial overdispersion requires point processes models
with finite dimensional distributions that are overdisperse relative to the
Poisson. Fitting such models usually heavily relies on the properties of
stationarity, ergodicity, and orderliness. And, though processes based on
negative binomial finite dimensional distributions have been widely
considered, they typically fail to simultaneously satisfy the three required
properties for fitting. Indeed, it has been conjectured by Diggle \& Milne
that no negative binomial model can satisfy all three properties. In light
of this, we change perspective, and construct a new process based on a
different overdisperse count model, the Generalized Waring Distribution.
While comparably tractable and flexible to negative binomial processes, the
Generalized Waring process is shown to possess all required properties, and
additionally span the negative binomial and Poisson processes as limiting
cases. In this sense, the GW process provides an approximate resolution to
the conundrum highlighted by Diggle \& Milne.

\bigskip

\textbf{Keywords and Phrases}: additivity; stationarity; ergodicity;
orderliness; overdispersion; Poisson process; negative binomial process;
generalized Waring process; complete separable metric space

\bigskip

\textbf{Running Head}: Spatial Overdispersion and the GWP

\section{Introduction}

\bigskip The definition of an appropriate probability model for spatial
count data typically requires the determination of an additive point process
on the domain in question. Additivity is a minimal requirement, requiring
that when the region of observation changes, or when non-overlapping regions
are aggregated in a systematic manner, the corresponding count distribution
remains in the same family. Additional assumptions that are often made for
convenience include stationarity/ergodicity (allowing estimation of the
model based on a single realization) and orderliness (to avoid the
apparition of coincident events).

In cases where the spatial counts to be modelled are prone to exhibit
overdispersion, however, it may be challenging to specify a point process
model that simultaneously features additivity, stationarity/ergodicity, and
orderliness. A popular modeling strategy is to construct a point process
with finite-dimensional laws of the negative binomial form. Such processes
are known as negative binomial processes, and have been defined and studied
on general state spaces (Gregoire (1983)). Owing to their combination of
flexibility and mathematical tractability, they have been employed in many
practical situations (see for example Bates (1955) , Boswell \& Patil (1977)
, Cliff\ \& Ord (1973), Ramakrishnan (1951) etc.). However, they have been
shown to fail in simultaneously accommodating the three properties listed
above. As a matter of fact, it has been conjectured by Diggle \& Milne
(1983), that additive/stationary/orderly spatial point processes processes
with negative binomial finite-dimensional distributions may not even exist.
In their words, it would seem that one is "unable to exhibit a negative
binomial point process that is statistically interesting according to the
criteria we laid down" [these criteria being
additivity/stationarity/orderliness].

To elaborate, the construction of a negative binomial process $N$ usually
hinges on one of two schemes. The first scheme is based on compounding
Poisson processes by means of the logarithmic distribution (see Feller
(1968). One defines $N(B)=\sum\limits_{k=1}^{M(B)}X_{k}$ to denote the count
corresponding to where $M$ is a stationary Poisson process with mean
(intensity) measure $E\left( N(B)\right) =\lambda \cdot \mu \left( B\right) $%
, \ where $\mu \left( B\right) $ denotes the area (Lebesgue measure) of $B$.
Given $M$, the random variables $X_{i}$ are taken to be independently and
identically distributed (i.i.d) according to the logarithmic series
distribution with parameters $\left( \delta ,z\right) $, having probability
generating function (p.g.f.) $\dfrac{-\ln (1-\frac{\delta z}{1+\delta })}{%
\ln (1+\delta )},\ \delta >0.$ The resulting process can be seen to be of
Negative binomial form with p.g.f. $E\{z^{N(B)}\}=\{1+\delta (1-z)\}^{\dfrac{%
-\lambda \cdot \mu \left( B\right) }{\ln (1+\delta )}}$. This is a Poisson
cluster scheme (cf. Daley \& Vere-Jones, 1972 , Example 2.4.B and Cox \&
Isham, 1980, Fisher, 1972\ , Example 5.6, Burnett \& Wasan, 1980), and, as
remarked by Diggle \& Milne,1983, is always stationary/ergodic (any
stationary Poisson cluster process is known to be mixing; Westcott, 1971, p.
300 ), but clearly non-orderly. A second scheme is based on mixing Poisson
processes, generating so-called Polya processes (see Matern (1971) , cf.
Daley \& Vere-Jones, (1972), Example 2.1.C; Fisher, (1972), p. 500). Here,
one samples a gamma random variable $\Lambda $ with parameters $\alpha $ and 
$\beta >0$, and conditionally specifies $N(B)$ to be Poisson given $\Lambda $%
, with intensity $\Lambda \cdot \mu \left( B\right) $. The resulting process
is again of the negative binomial type, with p.g.f. $E\{z^{N(B)}\}=\{1+\beta
(1-z)\}^{-\alpha \cdot \mu \left( B\right) }$. Polya processes on the real
line are well-established in the literature on accident proneness cf. Cane
(1972) . As mentioned again by Diggle \& Milne (1983), they are stationary
by construction. However, the only stationary mixed Poisson processes which
are ergodic are those for which the mixing distribution is concentrated at a
single point, thus giving an (ordinary) Poisson process (cf. Westcott, 1972,
p. 464). It follows that non-trivial processes of this type can be orderly
but never ergodic.

\bigskip In summary, the first approach yields ergodic but non-orderly
processes, whereas the second approach yields orderly but non-ergodic
processes. In this paper, therefore, rather than make a new attempt at
finding a point process with precisely negative binomial one-dimensional
distributions (which may not even be possible), we change strategy, and
consider a different choice of over-disperse one-dimensional distributions.
An established competitor to the negative binomial distribution is the
Generalized Waring Distribution (GWD; see, e.g. Irwin (1975), Xekalaki
(1983b, 1984)). This has long been used to fit overdisperse count data,
particularly in the field of accident studies, providing a more plausible
model for the interpretation of the data generating mechanism; and, it can
approximate the negative binomial and the Poisson distribution as limiting
cases. A corresponding (temporal) stochastic process has been defined and
studied by Xekalaki \& Zografi (2008), in the context of temporally evolving
data featuring clustering or contagion.

Using the GWD as a building block, we construct an additive, stationary,
ergodic, and orderly spatial point process, and study its basic properties.
We develop our results on a general separable metric state space, before
focussing on the practically relevant case of $\mathbb{R}^{d}$. The process
is seen to satisfy several useful closure properties (under projection,
marginalization, and superposition) and to be easy to simulate. We further
show that, in the limit as certain parameters of the process diverge, this
Generalized Waring Point Process approximates a negative binomial process.
In doing so, we give an approximate positive solution to the task set out by
Diggle \& Milne: while a stationary, ergodic and orderly \ point process
with one-dimensional negative binomial distributions may not exist, there
exists a point process that is stationary, ergodic and orderly point process
and has one dimensional distributions that are approximately negative
binomial (depending on parameter choice).

The paper is organised as follows. In Section 2, we provide some necessary
background notions related to the generalized Waring distribution, its
moments, and properties that will be used in subsequent sections.
Specifically, it is shown that the generalized Waring distribution posesses
the property of countable additivity, which is fundamental to our later
construction. The definition and existence of the generalized Waring process
in a complete separable metric space is given in Section 3. In particular,
the process is shown to be orderly, and to be characterised by the property
that $N\left( A\right) \ $ follows a Univariate Generalized Waring
Distribution (UGWD) with parameters $(a,k\mu \left( A\right) ,\rho )$ for
all bounded sets $A$ in a dissecting ring $\mathcal{A}$ of the complete
separable metric space. The same section includes the determination of the
corresponding intensity measure, factorial moment measures and the $n^{th}$
- order moment measures. The generalized Waring process in $%
\mathbb{R}
^{d}$ with Lebesgue measure as parameter measure $\mu \left( \cdot \right) $
is then defined in paragraph 4. It is shown to be orderly, ergodic and $n$%
th-order stationary. The existence of the $n$th-order reduced moments of a
generalized Waring process in $%
\mathbb{R}
^{d},$ if $\rho >n,$ useful for applications, is obtained as a corollary.
Finally, multivariate extensions are considered in Section 5, where we
define the multivariate GWP as a special case of the GWP on the product
space $S\times \left\{ 1,2,...m\right\} $, and is shown to satisfy several
appealing closure properties with respect to marginalization.

\section{\hspace{5mm}The Generalized Waring Distribution and Additivity}

In this section we provide some background on the generalized Waring
distribution and discuss some of its structural properties that will be
essential in what follows. In particular, we extend the previously
established finite additivity property to countable additivity, as a first
important step in the construction of the GW point process.

A random variable $X$ is said to have the generalized Waring distribution
with parameters $a,k$ and $\rho $, denoted by $GWD\left( a,k;\rho \right) $,
if 
\begin{equation}
P\left\{ X=n\right\} =\pi _{n}\left( a,k;\rho \right) =\frac{\rho _{\left(
k\right) }}{\left( \rho +a\right) _{\left( k\right) }}\frac{a_{\left(
n\right) }k_{\left( n\right) }}{\left( \rho +a+k\right) _{\left( n\right) }}%
\frac{1}{n!}\ \ \ _{n=0,1,...}  \label{1}
\end{equation}%
where $a_{\left( \beta \right) }=\dfrac{\Gamma \left( \alpha +\beta \right) 
}{\Gamma \left( \alpha \right) },$ $P\left( X=x\right) =0,$ $x\in
\{0,1,2,...\}^{C}$ (see e.g. Irwin (1975), Xekalaki (1981), \ Xekalaki
(1983b)). Here $a>0,$ $k>0,$ $\rho >0$ and $k$ need not be integers. The
distribution is symmetric in $a$ and $k$.

The probability generating function of the generalized Waring distribution
is given by 
\begin{equation}
E\left( z^{X}\right) =\overset{\infty }{\underset{n=0}{\sum }}z^{n}\pi
_{n}\left( a,k;\rho \right) =\frac{\rho _{\left( k\right) }}{\left( \rho
+a\right) _{\left( k\right) }}\ _{2}F_{1}(a,k;\rho +a+k;z)  \label{2}
\end{equation}%
\bigskip\ where%
\begin{equation*}
_{2}F_{1}(a,\beta ;\gamma ;z)=\overset{\infty }{\underset{n=0}{\sum }}\frac{%
a_{\left( n\right) }\beta _{\left( n\right) }}{\gamma _{\left( n\right) }}%
\frac{z^{n}}{n!}.
\end{equation*}%
The $r$th factorial moments are 
\begin{equation}
\mu _{\left[ r\right] }=\frac{a_{\left[ r\right] }k_{\left[ r\right] }}{%
\left( \rho -1\right) \left( \rho -2\right) ...\left( \rho -r\right) }
\label{3}
\end{equation}%
where $x_{\left[ r\right] }=x\left( x+1\right) ...\left( x+r\right) .$

From (\ref{3}) it follows immediately that all $r$th moments, ordinary
moments about any origin, central moments as well as factorial moments are
infinite if $\rho \leq r$. Moments about any origin, including central
moments, can be obtained from (\ref{3}) by the usual transformation formula
(see Irwin (1975), PartI ). In particular the mean is given by 
\begin{equation}
E\left( X\right) =\frac{ak}{\rho -1},\text{ }\rho >1  \label{4}
\end{equation}%
while the variance is 
\begin{equation}
\sigma ^{2}=\mu _{2}=\frac{ka\left( \rho +a-1\right) \left( \rho +k-1\right) 
}{\left( \rho -1\right) ^{2}\left( \rho -2\right) },\text{ }\rho >2
\label{5}
\end{equation}

The multivariate generalized Waring distribution with parameter vector $%
(\alpha ,$ $k_{1},$ $...,$ $k_{s};$ $\rho )$, denoted by $MGWD$ $(a;$ $%
\mathbf{k};$ $\rho )$, is the probability distribution of a random vector $%
\left( X_{i},\ i=1,\text{ }2,\text{ }...,\text{ }s\right) $ of nonnegative
integer-valued components, with probability function given by 
\begin{equation}
P_{x_{1},...,x_{s}}=P\left( X_{i}=x_{i},i=1,2,...,s\right) =\frac{\rho
_{\left( \overset{s}{\underset{i=1}{\sum }}k_{i}\right) }a_{\left( \overset{s%
}{\underset{i=1}{\sum }}x_{i}\right) }}{\left( \rho +a\right) _{\left( 
\overset{s}{\underset{i=1}{\sum }}k_{i}+\overset{s}{\underset{i=1}{\sum }}%
x_{i}\right) }}\overset{s}{\underset{i=1}{\tprod }}\frac{k_{i\left(
x_{i}\right) }}{x_{i}!}  \label{6}
\end{equation}%
(see Xekalaki (1986)). The special case for $s=2$ is known in the literature
as the bivariate generalized Waring distribution, denoted by $BGWD\left(
a;k_{1},k_{2};\rho \right) .$

\bigskip The probability generating function of the multivariate Generalized
Waring distribution can be expressed in terms of Lauricella's hypergeometric
function of type D as 
\begin{equation*}
G\left( \underline{z}\right) =\frac{\rho _{\left( \overset{s}{\underset{i=1}{%
\sum }}k_{i}\right) }}{\left( \rho +a\right) _{\left( \overset{s}{\underset{%
i=1}{\sum }}k_{i}\right) }}\ F_{D}(a,k_{1},k_{2},...,k_{s};\rho +a+\overset{s%
}{\underset{i=1}{\sum }}k_{i};\underline{z})
\end{equation*}%
where $\ $%
\begin{equation*}
F_{D}(a,\beta _{1},\beta _{2},...,\beta _{s};\gamma ;\underline{z})=\underset%
{r_{1},r_{2},...r_{s}}{\sum }\frac{a_{\left( \tsum r_{i}\right) }}{\gamma
_{\left( \tsum r_{i}\right) }}\overset{s}{\underset{i=1}{\tprod }}\frac{%
\left( \beta _{i}\right) _{\left( r_{i}\right) }\text{\ }\left( z_{i}\right)
^{r_{i}}}{r_{i}!},
\end{equation*}

\begin{equation*}
\underline{z}=(z_{1},\text{ }z_{2},\text{ }...,\text{ }z_{s})
\end{equation*}

$\bigskip $The factorial moments of the $MGWD\left( a;\text{ }\mathbf{k};%
\text{ }\rho \right) $ (see Xekalaki (1985a), Xekalaki (1986)) are then
given by%
\begin{eqnarray}
\mu _{\left( r_{1},r_{2},...r_{l}\right) } &=&E\left[ \left( X_{1}\right) _{%
\left[ r_{1}\right] }\left( X_{2}\right) _{\left[ r_{2}\right] }...\left(
X_{s}\right) _{\left[ r_{s}\right] }\right]  \label{7} \\
&=&\frac{a_{\left( \tsum r_{i}\right) }\overset{s}{\underset{i=1}{\tprod }}%
\left( k_{i}\right) _{\left( r_{i}\right) }}{\left( \rho -1\right) \left(
\rho -2\right) ...\left( \rho -\tsum r_{i}\right) },\text{ }r_{i}=0,1,...;%
\text{ }i=1,2,...,s
\end{eqnarray}%
and are finite for $\rho >\tsum r_{i}$, the latter being a necessary
condition for the series $F_{D}(a,k_{1}+r_{1},k_{2}+r_{2},...,k_{s}+r_{s};%
\rho +a+\overset{s}{\underset{i=1}{\sum }}\left( k_{i}+r_{i}\right) ;%
\underline{1})$ to converge. Moments of order $n$ can be derived from these
factorial moments.

The marginal means and marginal variances are respectively given by 
\begin{equation}
\mu _{X_{i}}=E\left( X_{i}\right) =\frac{ak_{i}}{\rho -1},\text{ }\rho >1
\label{8}
\end{equation}%
$\bigskip $ 
\begin{equation}
\sigma _{X_{i}}^{2}=\frac{k_{i}a\left( \rho +a-1\right) \left( \rho
+k_{i}-1\right) }{\left( \rho -1\right) ^{2}\left( \rho -2\right) },\rho >2
\label{9}
\end{equation}%
$i=1,2,...,s.$(see Xekalaki (1986)

The second moment and the pairwise covariances are 
\begin{equation}
\mu _{X_{i}X_{j}}=E\left( X_{i}X_{j}\right) =\frac{a\left( a+1\right)
k_{i}k_{j}}{\left( \rho -1\right) \left( \rho -2\right) },\text{ }%
i,j=1,2,...,s;\text{ }\rho >2  \label{10}
\end{equation}

\begin{equation}
\sigma _{X_{i}X_{j}}=\frac{a\left( \rho +a-1\right) k_{i}k_{j}}{\left( \rho
-1\right) ^{2}\left( \rho -2\right) },i,j=1,2,...,s;\rho >2  \label{11}
\end{equation}

One of the most important features of the GWD is \textbf{additivity.}
Specifically, if $X$ and $Y$ are random variables with marginal
distributions $UGWD$ $(a,$ $k_{1};$ $\rho )$ and $UGWD$ $(a,$ $k_{2};$ $\rho
)$, respectively, and with joint distribution $BGWD$ $(a;$ $k_{1},$ $k_{2};$ 
$\rho )$, then $X+Y$ is a $UGWD$ $(a,$ $k_{1}+k_{2};$ $\rho )$ random
variable. More generally, letting $X_{j}$ be $UGWD$ $(a,$ $k_{j};$ $\rho )$
or each $j,\ j=1,2,...,n\ $ \ and jointly distributed as $MGWD$ $(a;$ $%
k_{1}, $ $k_{2},$ $...,$ $k_{n};$ $\rho )$, then, if we denote $m=\overset{n}%
{\underset{j=1}{\sum }}k_{j},$ we have that $S=\overset{n}{\underset{j=1}{%
\sum }}X_{j}$ also has a $UGWD$ $(a,$ $m;$ $\rho )$ distribution.

These last two properties hint at the possibility of using the GWD as a
basis for the construction of overdisperse point processes. This requires
extending additivity to countable additivity, which we do in the form of the
next theorem:

\begin{theorem}
\label{cath}Let $X_{j}$ be $UGWD\left( a,k_{j};\rho \right) $ variables for
each $j,\ j=1,2,...\ $ \ and for each $n\geq 3$ let their joint distribution
be the $MGWD\left( a,k_{1},k_{2},...,k_{n};\rho \right) .$ If $m=\overset{%
\infty }{\underset{j=1}{\sum }}k_{j}$ converges, then $S=\overset{\infty }{%
\underset{j=1}{\sum }}X_{j}$ converges with probability $1$, and $S$ has a $%
UGWD\left( a,m;\rho \right) $ distribution. If on the other hand, $\overset{%
\infty }{\underset{j=1}{\sum }}$ $k_{j}$ diverges, then $S$ diverges with
probability $1$.\bigskip
\end{theorem}

\textit{Proof} By induction on $n,$ the random variable \ $S_{n}=\overset{n}{%
\underset{j=0}{\sum }}X_{j}$ has a $UGWD$ $(a,$ $m_{n};$ $\rho )$
distribution, where $m_{n}=\overset{n}{\underset{j=1}{\sum }}$ $k_{j}.$Thus,
for any $r$, 
\begin{equation}
P\left\{ S_{n}\leq r\right\} =\overset{r}{\underset{i=0}{\sum }}\pi
_{i}\left( a,m_{n};\rho \right)  \label{12}
\end{equation}%
The sequence $\left\{ S_{n}\leq r\right\} $ is a decreasing sequence of
events for fixed $r$, and their intersection is $\left\{ S\leq r\right\} .$%
Thus, using continuity from above, 
\begin{eqnarray*}
P\left\{ S\leq r\right\} &=&\underset{n\rightarrow \infty }{\lim }P\left\{
S_{n}\leq r\right\} \\
&=&\underset{n\rightarrow \infty }{\lim }\overset{r}{\underset{i=0}{\sum }}%
\pi _{i}\left( a,m_{n};\rho \right) .
\end{eqnarray*}%
If $m_{n}$ converges to a finite limit $m,$ the continuity of $\pi _{j}$
implies that 
\begin{equation}
P\left\{ S\leq r\right\} =\overset{r}{\underset{i=0}{\sum }}\pi _{i}\left(
a,m;\rho \right)  \label{13}
\end{equation}%
leading to 
\begin{equation}
P\left\{ S=r\right\} =\pi _{r}\left( a,m;\rho \right) .  \label{14}
\end{equation}%
This in turn implies that $S$ is finite and distributed as generalized
Waring with parametrs $a,m;\rho $ ($UGWD\left( a,m;\rho \right) $).

On the other hand, if $m_{n}\rightarrow \infty ,$\qquad \qquad \qquad 
\begin{equation*}
\overset{r}{\underset{i=0}{\sum }}\pi _{i}\left( a,m_{n};\rho \right) =%
\overset{r}{\underset{i=0}{\sum }}\frac{\rho _{\left( m_{n}\right) }}{\left(
\rho +a\right) _{\left( m_{n}\right) }}\frac{a_{\left( i\right) }m_{n\left(
i\right) }}{\left( \rho +a+m_{n}\right) _{\left( i\right) }}\frac{1}{i!}=
\end{equation*}%
\begin{equation*}
\left[ \frac{\rho _{\left( a\right) }}{\left( \rho +m_{n}\right) _{\left(
a\right) }}\frac{a_{\left( i\right) }}{i!}\frac{m_{n}}{\left( \rho
+a+m_{n}\right) }\frac{\left( m_{n}+1\right) }{\left( \rho +a+m_{n}+1\right) 
}...\frac{\left( m_{n}+i-1\right) }{\left( \rho +a+m_{n}+i-1\right) }%
\rightarrow 0\text{ }\right]
\end{equation*}%
so that $P\left\{ S>r\right\} =1.$ Since this holds for all $r$, $S$
diverges with probability $1.\bigskip $

\section{The Generalized Waring Process}

We now proceed to the definition of the generalized Waring process on a
complete separable metric space and the investigation of some of its basic
properties. The construction starts from postulating the existence of a
point process with finite dimensional distributions of the generalized
Waring form (Subsection 3.1), and then demonstrating the existence and
uniqueness of such a process (Subsection 3.2). Basic features of the process
such as a conditional property useful for simulation, as well as its
intensity measure, factorial moment measures and $n$th - order moment
measures are then derived in Subsection 3.3.

\subsection{Definition and Basic Properties}

Let $\mathcal{S}$ be a complete separable metric space, $\mathcal{A}$ a
semiring of bounded Borel sets generating the Borel $\sigma $-algebra $%
\mathcal{B}_{S}$ of subsets of $\mathcal{S}$ (Appendix2. Lemma A2.I.III,
Daley and Vere-Jones (1988)) and $\mu \left( \cdot \right) $ a boundedly
finite Borel measure. The distribution of a random measure is completely
determined by its finite dimensional (fidi) distributions, i.e. the joint
distribution of arbitrary finite families $\{A_{i},$ $i=1,$ $...,$ $s\}$ of
disjont sets from $\mathcal{A}$ (Proposition 6.2.III, Daley \& Vere-Jones
(1988)). Now consider the space of all boundedly finite, integer-valued
measures $(\overset{\wedge }{\mathcal{N}_{S}},$ $\mathcal{B}\left( \overset{%
\wedge }{\mathcal{N}_{S}}\right) )$ and let $\left( \Omega ,\mathcal{F},%
\mathcal{P}\right) $ be some probability space.

\textbf{Definition 1} Let 
\begin{equation*}
N:\left( \Omega ,\mathcal{F},\mathcal{P}\right) \rightarrow \left( \overset{%
\wedge }{\mathcal{N}_{S}},\mathcal{B}\left( \overset{\wedge }{\mathcal{N}_{S}%
}\right) \right)
\end{equation*}%
be a point process for whose finite dimensional distributions over disjoint
bounded Borel sets $\{A_{i},$ $i=1,$ $...,$ $l$ $\}$ \ are given by 
\begin{equation}
P\left\{ N\left( A_{i}\right) =n_{i};i=1,...,l\right\} =\frac{\rho _{\left( k%
\overset{l}{\underset{i=1}{\sum }}\mu \left( A_{i}\right) \right) }a_{\left( 
\overset{l}{\underset{i=1}{\sum }}n_{i}\right) }}{\left( \rho +a\right)
_{\left( k\overset{l}{\underset{i=1}{\sum }}\mu \left( A_{i}\right) +\overset%
{l}{\underset{i=1}{\sum }}n_{i}\right) }}\overset{l}{\underset{i=1}{\tprod }}%
\frac{\left[ k\mu \left( A_{i}\right) \right] _{\left( n_{i}\right) }}{n_{i}!%
}.  \label{1a}
\end{equation}%
Then $N$ is called a \textit{generalized Waring process} with parameters $a$%
, $\rho ,$ $k>0\ $and\ parameter measure $\mu \left( \cdot \right) $.\bigskip

In other words, for every finite family of disjoint bounded Borel sets\ $%
\{A_{i},$ $i=1,$ $...,$ $l$\bigskip $\}$ the joint distribution of $\left\{
N\left( A_{i}\right) =n_{i},i=1,...,l\right\} $ is the $MGWD$ $(a,$ $k\mu
\left( A_{1}\right) ,$ $k\mu \left( A_{2}\right) ,$ $...,$ $k\mu \left(
A_{l}\right) ;$ $\rho ).$ As usual, the process $\left\{ N\left( A\right) ;%
\text{ }A\in \mathcal{B}_{S}\right\} $ is to be thought of as a random
measure. In particular, for any $A\in \mathcal{B}_{S}$, $N\left( A\right) $
is a $Z^{t}-$valued random variable, while for any $\omega \in \Omega ,$ $%
N\left( \omega ,\cdot \right) $ is a discrete Radon measure.

We remark that, if such a process exists, it will necessarily be countably
additive. To see this, let $\left\{ A_{i},i=1,2,...\right\} $ be disjoint
and have union $A$. Using Theorem 1, and the fact that $\mu \left( A\right) =%
\overset{\infty }{\underset{i=1}{\sum }}\mu \left( A_{i}\right) $ converges,
we immediately obtain that $N\left( A\right) =\overset{\infty }{\underset{i=1%
}{\sum }}N\left( A_{i}\right) $ is distributed as $UGWD\left( a,k\mu \left(
A\right) ;\rho \right) $. Furthermore, such a process will be orderly
provided the parameter measure is diffuse:

\begin{theorem}
\textit{\ A process as in Definition 1 is an orderly point process if and
only if its parameter measure has no fixed atoms.}
\end{theorem}

\textit{Proof} \bigskip A point process is orderly when given any bounded $%
A\in \mathcal{B}_{S},$ there is a dissecting system $\mathcal{T=}\left\{ 
\mathcal{T}_{n}\right\} =\{$ $\{$ $A_{ni}:$ $i=1,$ $...,$ $k_{n}$ $\}$ $\}$
such that $\underset{\mathcal{T}_{n}}{\inf }$ \ $\sum_{i-1}^{k_{n}}$ $%
P\left\{ N\left( A_{ni}\right) >2\right\} =0.$ (see Daley \& Vere-Jones
(1988)). Hence it is sufficient to examine when the ratio $P$ $\{N\left(
A_{\varepsilon ,x}\right) $ $>$ $1\}/$ $P$ $\{N\left( A_{\varepsilon
,x}\right) $ $>$ $0\}$ tends to $0,$ where $A_{\varepsilon ,x}$ is the open
sphere of radius $\varepsilon $ and center $x\in A$. In the case of a GW
process, $N\left( A_{\varepsilon ,x}\right) $ has a generalized Waring
distribution with parameters $a>0$, $\rho >0$ and $\mu \left( A_{\varepsilon
,x}\right) =\mu _{\varepsilon }$, so that 
\begin{equation*}
P\left\{ N\left( A_{\varepsilon ,x}\right) >0\right\} =1-P\left\{ N\left(
A_{\varepsilon ,x}\right) =0\right\} =1-\frac{\rho _{\left( k\mu
_{\varepsilon }\right) }}{\left( \rho +a\right) _{\left( k\mu _{\varepsilon
}\right) }},
\end{equation*}%
\begin{equation*}
P\left\{ N\left( A_{\varepsilon ,x}\right) >1\right\} =1-\frac{\rho _{\left(
k\mu _{\varepsilon }\right) }}{\left( \rho +a\right) _{\left( k\mu
_{\varepsilon }\right) }}-\frac{\rho _{\left( k\mu _{\varepsilon }\right) }}{%
\left( \rho +a\right) _{\left( k\mu _{\varepsilon }\right) }}\frac{a\cdot
k\mu _{\varepsilon }}{\left( \rho +a+k\mu _{\varepsilon }\right) }.
\end{equation*}%
If $x$ is a fixed atom of $\mu $, then $\mu _{\varepsilon }\rightarrow \mu
_{0}=\mu \left\{ x\right\} >0$ as $\varepsilon \rightarrow 0,$ while if $x$
is not a fixed atom, then $\mu \left( A_{\varepsilon ,x}\right) \rightarrow
0.$

In the first case, the ratio $P\left\{ N\left( A_{\varepsilon ,x}\right)
>1\right\} /P\left\{ N\left( A_{\varepsilon ,x}\right) >0\right\} $ tends to
the constant $1-\dfrac{\rho _{\left( k\mu _{0}\right) }\cdot a\cdot k\mu _{0}%
}{\left( \rho +a\right) _{\left( k\mu _{0}+1\right) }-\rho _{\left( k\mu
_{0}\right) }},$ while in the second case it tends to $0$, and the proof is
complete.

From now and on we will consider only orderly generalized Waring processes.
Indeed, any orderly point process with finite dimensional distributions of
the generalized Waring type is necessarily a GWP with a non-atomic parameter
measure:

\begin{theorem}
\label{nsc}Let $N\left( \cdot \right) $ be an orderly point process. For $%
N\left( \cdot \right) $ to be a generalized Waring process with parameters $%
a>0$, $\rho >0,$ $k>0$ and parameter measure $\mu \left( \cdot \right) $, it
is necessary and sufficient that there exist a boundedly finite nonatomic
measure $\mu $ on the Borel sets $B_{s}$ such that $N\left( A\right) $ has
generalized Waring distribution with parameters $a,k\mu \left( A\right)
,\rho $ for each bounded set $A\ $of a dissecting ring $\mathcal{A}$ of the
complete separable metric space $\mathcal{S}$.
\end{theorem}

\textit{Proof}\ We begin with necessity. Let $N\left( \cdot \right) $ be a
generalized Waring Process and $A$ a\ bounded set of a dissecting ring $%
\mathcal{A}$ ($A$ is also a Borel set). Then, by definition, there exists a
boundedly finite Borel measure $\mu \left( \cdot \right) $ such that for
every finite family of disjoint bounded Borel sets $\{A_{i},$ $i=1,...,s\},$ 
$P\{N\left( A_{i}\right) =n_{i},$ $i=1,...,s\}$ is given by \ref{1a}. From
this, it follows that the distribution of $N\left( A\right) $ is the $%
GWD\left( a,\text{ }k\mu \left( A\right) ;\text{ }\rho \right) .$

To prove sufficiency, suppose that there exists a boundedly finite nonatomic
measure $\mu $ on the Borel sets $B_{s}$ such that $N\left( A\right) $ has
generalized Waring distribution with parameter $a,k\mu \left( A\right) ,\rho 
$ for each bounded set $A$ of a dissecting ring. According to Theorem 7.3.II
of Daley \& Vere-Jones (1988), the values of the avoidance function $%
P_{0}\left( A\right) =$ $P\left\{ N\left( A\right) =0\right\} $ $=\dfrac{%
\rho _{\left( k\mu \left( A\right) \right) }}{\left( \rho +a\right) _{\left(
k\mu \left( A\right) \right) }}$ \ on the bounded sets of a dissecting ring
for$\ $the complete separable metric space, determine the distribution of a
simple point process $N\left( \cdot \right) $ on this space.

\subsection{\protect\bigskip Existence and Uniqueness}

To prove that the point process stipulated in the previous section does
indeed exist, it is sufficient to establish that the fidi distributions
given by (\ref{1a}) fulfill Kolmogorov's consistency conditions, combined
with the measure requirements given by the basic existence theorem for point
processes (Theorem 7.I.XI Daley \& Vere-Jones (1988)).

\begin{theorem}
(\textit{Kolmogorov's Consistency Conditions) } A collection of finite
dimensional distributions as defined via Definition 2 satisfies Kolmogorov's
consistency conditions. That is, for every finite family of disjoint bounded
Borel sets $\left\{ A_{i},i=1,...,l\right\}$,

\bigskip (I) for any permutation $i_{1},...,i_{l}$ of the indexes $1,...,l$ 
\begin{equation}
P_{l}\left( A_{i_{1}},...,A_{i_{l}};n_{i_{1}},...,n_{i_{l}}\right)
=P_{l}\left( A_{1},...,A_{l};n_{1},...,n_{l}\right)  \label{2a}
\end{equation}

\bigskip (II) $\underset{r=0}{\overset{\infty }{\sum }}P_{l}\left(
A_{1},...,A_{l},n_{1},...,n_{l-1},r\right) =P_{l-1}\left(
A_{1},...,A_{l-1},n_{1},...,n_{l-1}\right) $\bigskip
\end{theorem}

\textit{Proof} To show (I), we notice that one can write $\overset{l}{%
\underset{j=1}{\sum }}\mu $ $=$ $\overset{l}{\underset{j=1}{\sum }}\mu
\left( A_{j}\right) ,$ $\overset{l}{\underset{j=1}{\sum }}n_{i_{j}}$ $=%
\overset{l}{\underset{j=1}{\sum }}n_{j},$ $\overset{l}{\underset{j=1}{\tprod 
}}$ $\dfrac{\left[ k\mu \left( A_{i_{j}}\right) \right] _{\left(
n_{i_{j}}\right) }}{n_{i}!}$ $=\overset{l}{\underset{j=1}{\tprod }}\dfrac{%
\left[ k\mu \left( A_{j}\right) \right] _{\left( n_{j}\right) }}{n_{j}!}$
which proves (\ref{2a}).

To show (II), we write $\underset{r=0}{\overset{\infty }{\sum }}$ $P_{l}($ $%
A_{1},$ $...,$ $A_{l};$ $n_{1},$ $...,$ $n_{l-1},$ $r$ $)$ $=$ $P_{l-1}($ $%
A_{1},$ $...,$ $A_{l-1};$ $n_{1},$ $...,$ $n_{l-1})$ $\underset{r=0}{\overset%
{\infty }{\sum }}$ $\dfrac{\left( \rho +k\overset{l-1}{\underset{i=1}{\sum }}%
\mu \left( A_{i}\right) \right) _{\left( k\mu \left( A_{r}\right) \right)
}\left( a+\overset{l-1}{\underset{i=1}{\sum }}n_{i}\right) _{\left( r\right)
}}{\left( \rho +a+k\overset{l-1}{\underset{i=1}{\sum }}\mu \left(
A_{i}\right) +\overset{l-1}{\underset{i=1}{\sum }}n_{i}\right) _{\left( k\mu
\left( A_{r}\right) +r\right) }}$ $\dfrac{\left[ k\mu \left( A_{i}\right) %
\right] _{\left( r\right) }}{r!}$ $=P_{l-1}($ $A_{1},$ $...,$ $%
A_{l-1};n_{1}, $ $...,$ $n_{l-1})$

\begin{theorem}
\bigskip \label{ex}\bigskip (\textit{Measure Requirements}) Suppose that

\bigskip $(I)$ $N$ is bounded finite a.s. and has no fixed atoms.

\bigskip (II) $N$ satisfies Definition 2.

\bigskip Then, there exists a boundedly finite nonatomic Borel measure $\mu
\left( \cdot \right) $ such that $P_{0}\left( A\right) =\Pr \left\{ N\left(
A\right) =0\right\} =\dfrac{\rho _{\left( k\mu \left( A\right) \right) }}{%
\left( \rho +a\right) _{\left( k\mu \left( A\right) \right) }}$ for all
bounded borel sets $A$ and $\forall i,i=1,...,s$ $\mu \left( A_{i}\right)
=\mu _{i}$.\bigskip
\end{theorem}

\textit{Proof} \bigskip Let $A\in B_{s}\ $and let $\mu \left( A\right) >0$
be the root of the equation $P_{0}\left( A\right) =\dfrac{\rho _{\left( k\mu
\left( A\right) \right) }}{\left( \rho +a\right) _{\left( k\mu \left(
A\right) \right) }}$ which does exist (see Appendix, Lemma \ref{exl}).

a) We first prove that $\mu \left( \cdot \right) $ is a measure. To show
finite additivity, we observe that 
\begin{equation*}
P_{0}\left( A\right) =\Pr \left\{ N\left( A\right) =0\right\} =\frac{\rho
_{\left( k\mu \left( A\right) \right) }}{\left( \rho +a\right) _{\left( k\mu
\left( A\right) \right) }}.
\end{equation*}%
Hence for each family of bounded, disjoint, Borel sets $\left\{
A_{i},i=1,...,s\right\} $, the joint distribution of $\left\{ N\left(
A_{i}\right) =n_{i},i=1,...,s\right\} $ is the $MGWD($ $a,$ $k\mu \left(
A_{1}\right) ,$ $k\mu \left( A_{2}\right) ,$ $...,$ $k\mu \left(
A_{s}\right) ;$ $\rho $ $),$ and if $A=\overset{s}{\underset{i=1}{\sum }}%
A_{i}$ then $N\left( A\right) =\overset{s}{\underset{i=1}{\sum }}N\left(
A_{i}\right) $ has distribution $GWD\left( a,k\mu \left( A\right) ;\rho
\right) .$ So $\mu \left( A\right) =\overset{s}{\underset{i=1}{\sum }}\mu
\left( A_{i}\right) $ which establishes finite additivity of $\mu \left(
\cdot \right) .$ To extend this to countable additivity, it suffices to
prove that $\mu \left( A_{i}\right) \rightarrow 0$ for any decreasing
sequence $\left\{ A_{i}\right\} $ of bounded Borel sets for which $\mu
\left( A_{i}\right) <\infty $ and $A_{i}\downarrow $\O . For $%
A_{i}\downarrow $\O\ $N\left( A_{i}\right) \rightarrow 0$ a.s. and thus $%
P_{0}\left( A_{i}\right) =\Pr \left\{ N\left( A_{i}\right) =0\right\}
\rightarrow 1$ a.s. hence $\mu \left( A_{i}\right) =\dfrac{\rho \left(
1-P_{0}\left( A_{i}\right) \right) }{kP_{0}\left( A_{i}\right) }\rightarrow
0 $ a.s.

\bigskip b)To show that $\mu \left( \cdot \right) $ is non-atomic, we can
consider by $\left( I\right) $ that for every $x$ that $\Pr \left\{ \text{ }%
N\left( \left\{ x\right\} \right) >0\text{ }\right\} $ $\ =$ $\left(
1-P_{0}\left( \left\{ x\right\} \right) \right) =0.$So $\mu \left( \left\{
x\right\} \right) =\dfrac{\rho \left( 1-P_{0}\left( \left\{ x\right\}
\right) \right) }{kP_{0}\left( \left\{ x\right\} \right) }=0$\bigskip

c)To show that $\mu \left( \cdot \right) $ is boundedly finite it is enough
to prove that $P_{0}\left( A\right) >0$ for every bounded borel set $A.$ By
supposing the contrary that for some set $A$, $P_{0}\left( A\right) =0,$
one, following Daley \& Vere-Jones (1988),\ Lemma 2.4.VI, can find that in
this case there exists a fixed atom of the process, contradicting $\left(
I\right) $ which proves that $P_{0}\left( A\right) >0$ for every bounded
borel set $A.$\bigskip\ 

\subsection{Conditional Property and Moment Measures}

A useful property of the GWP is the conditional property, which provides a
straightforward way of simulating the process:

\begin{theorem}
(Conditional Property). Consider a Generalized Waring point process in $%
\Omega $ with parameters $a>0,\rho >0,k>0$ . Let $W\subset \Omega $ be any
region with \ $0<\mu \left( W\right) <+\infty $. Given that $N(W)=n$, the
conditional distribution of $N(B)$ for $B\subset W$ is the beta-binomial
distribution with parameters $\mu \left( B\right) ,$ $\mu \left( W)-\mu
(B\right) $ and $n$ :%
\begin{equation*}
p\left( N\left( B\right) =k\left\vert N\left( W\right) =n\right. \right) =%
\dbinom{n}{k}\dfrac{\left( \mu \left( B\right) \right) _{\left( k\right)
}\left( \mu \left( W)-\mu (B\right) \right) _{\left( n-k\right) }}{\left(
\mu \left( W\right) \right) _{\left( n\right) }}
\end{equation*}
\end{theorem}

$\bigskip $

\textit{Proof} $\bigskip $

$\qquad \qquad p\left( N\left( B\right) =k\left\vert N\left( W\right)
=n\right. \right) =\dfrac{p\left( N\left( B\right) =k,N\left( W-B\right)
=n-k\right) }{p\left( N\left( W\right) =n\right) }$

$\qquad =\dfrac{\dfrac{\rho _{\left( a\right) }}{\left( \rho +\mu \left(
W\right) \right) _{\left( a\right) }}\dfrac{a_{\left( n\right) }\left( \mu
\left( B\right) \right) _{\left( k\right) }\left( \mu \left( W-B\right)
\right) _{\left( n-k\right) }}{\left( \rho +\mu \left( W\right) +a\right)
_{\left( n\right) }}\dfrac{1}{k!}\dfrac{1}{\left( n-k\right) !}}{\dfrac{\rho
_{\left( a\right) }}{\left( \rho +\mu \left( W\right) \right) _{\left(
a\right) }}\dfrac{a_{\left( n\right) }\left( \mu \left( W\right) \right)
_{\left( n\right) }}{\left( \rho +\mu \left( W\right) +a\right) _{\left(
n\right) }}\dfrac{1}{n!}}$

$\qquad =\dfrac{n!}{k!\left( n-k\right) !}\dfrac{\left( \mu \left( B\right)
\right) _{\left( k\right) }\left( \mu \left( W)-\mu (B\right) \right)
_{\left( n-k\right) }}{\left( \mu \left( W\right) \right) _{\left( n\right) }%
}$

$\qquad =\dbinom{n}{k}\dfrac{\left( \mu \left( B\right) \right) _{\left(
k\right) }\left( \mu \left( W)-\mu (B\right) \right) _{\left( n-k\right) }}{%
\left( \mu \left( W\right) \right) _{\left( n\right) }}$

\bigskip Using the conditional property, we can generate a realization of a
Generalized Waring process with parameters $a>0,\rho >0,k>0$ in $W$, through
the following steps:

1. Generate a random variable $M$ with a Generalized Waring distribution
with parameters $a,\rho ,k\cdot \mu \left( W\right) .$

2. Given $M=m$, generate m points $Z_{1},Z_{2},...,Z_{m}$ in $W$ where $Z_{i}%
\symbol{126}Bernoulli(\pi _{i})$ having first simulated a draw from the beta
process i.e. a countably infinite collection of weighted atoms in $W$ , with
weights that lie in the interval [0; 1] (Hjort (1990)).

\bigskip We now turn our attention to determining the $n^{th}$ - order
moment measures of the process, needed to establish $n^{th}$ order
stationary, as discussed in the next section.

Let $N$ be a generalized Waring process with parameters $\left( a,k;\rho
\right) $ and parameter measure $\mu \left( \cdot \right) $. For $A$ a Borel
set, the distribution of $N\left( A\right) $ is the $GWD\left( a,k\mu \left(
A\right) ;\rho \right) .$ Therefore, its first moment measure is 
\begin{equation*}
\lambda \left( A\right) =E\left( N\left( A\right) \right) =\frac{ak\mu
\left( A\right) }{\rho -1},\text{ }\rho >1
\end{equation*}%
and its intensity rate is the Radon-Nikodym derivative 
\begin{equation*}
\eta \left( A\right) =\frac{d\lambda }{d\mu }=\frac{ak}{\rho -1},\text{ }%
\rho >1.
\end{equation*}%
For $A,$ $B$ two Borel sets the joint distribution of $\left( N\left(
A\right) ,N\left( B\right) \right) $ is the $BGWD($ $a,k\mu \left( A\right)
,k\mu \left( B\right) ;\rho $ $),$ hence the second-order moment measure of
the process is 
\begin{equation}
M_{2}\left( A\times B\right) =E\left( N\left( A\right) N\left( B\right)
\right) =\frac{a\left( a+1\right) k^{2}\mu \left( A\right) \mu \left(
B\right) }{\left( \rho -1\right) \left( \rho -2\right) },\text{ }\rho >2
\label{16}
\end{equation}

Given a finite family of disjoint bounded Borel sets $\left\{
A_{i},i=1,...,s\right\} $\bigskip\ the joint distribution of $\left\{
N\left( A_{i}\right) =n_{i},i=1,...,s\right\} $ is the $MGWD$ $($ $a,$ $k\mu
\left( A_{1}\right) ,$ $k\mu \left( A_{2}\right) ,$ $...,$ $k\mu \left(
A_{s}\right) $ $;\rho $ $)$, hence the factorial moment measure, $E$ $%
[N\left( A_{1}\right) _{\left[ r_{1}\right] }$ $N\left( A_{2}\right) _{\left[
r_{2}\right] }$ $...$ $N\left( A_{s}\right) _{\left[ r_{s}\right] }],$ of
the process is 
\begin{equation}
\mu _{\left( r_{1},r_{2},...r_{l}\right) }\left( A_{1}\times A_{2}\times
...\times A_{s}\right) =\dfrac{a_{\left( \tsum r_{i}\right) }k^{s}\overset{s}%
{\underset{i=1}{\tprod }}\left( \mu \left( A_{i}\right) \right) _{\left(
r_{i}\right) }}{\left( \rho -1\right) \left( \rho -2\right) ...\left( \rho
-\tsum r_{i}\right) },  \label{17q}
\end{equation}
$r_{i}=0,1,...;$ $i=1,2,...,s.$ The $n$th - order moment measures can now be
obtained from (\ref{17q}) 
\begin{eqnarray}
M_{n}\left( A_{1}\times A_{2}\times ...\times A_{n}\right) &=&E\left[ \left(
N\left( A_{1}\right) \right) _{\left[ 1\right] }\left( N\left( A_{2}\right)
\right) _{\left[ 1\right] }...\left( N\left( A_{s}\right) \right) _{\left[ 1%
\right] }\right]  \label{18} \\
&=&\frac{a_{\left( n\right) }k^{n}\overset{n}{\underset{i=1}{\tprod }}\left(
\mu \left( A_{i}\right) \right) _{\left( r_{i}\right) }}{\left( \rho
-1\right) \left( \rho -2\right) ...\left( \rho -n\right) },\text{ for }\rho
>n.  \notag
\end{eqnarray}%
for $\rho >\tsum r_{i}.$

\section{The Generalized Waring Process in $%
\mathbb{R}
^{d}$}

We now focus on the generalized Waring process on the state-space $%
\mathbb{R}
^{d}$, with Lebesgue measure as its parameter measure $\mu \left( \cdot
\right) $. We show that this constitutes an orderly, stationary, ergodic and 
$n$th-order stationary point process.

\subsection{The Generalized Waring Process as a Simple Point Process}

Let $S=%
\mathbb{R}
^{d}$ and let $\mu \left( \cdot \right) $ be the Lebesgue measure on $%
\mathbb{R}
^{d}$\bigskip .\ The Borel algebra $\mathcal{B}_{%
\mathbb{R}
^{d}}$ in $%
\mathbb{R}
^{d}$ is the smallest $\sigma -algebra$ on $%
\mathbb{R}
^{d}$ which contains all the open rectangles of $d-$dimensions. The
generalized Waring process $\{N\left( A\right) ;$ $A\in \mathcal{B}_{%
\mathbb{R}
^{d}}\}$ can be defined by assuming that for every finite family of disjoint
bounded Borel sets $\{A_{i},$ $i=1,$ $...,$ $s\}$ the joint distribution of $%
\{N\left( A_{i}\right) =$ $n_{i},$ $i=1,$ $...,$ $s\}$ is the $MGWD$ $(a,$ $%
k\mu \left( A_{1}\right) ,$ $k\mu \left( A_{2}\right) ,$ $...,$ $k\mu \left(
A_{s}\right) ;$ $\rho )$ $,a>0,$ $\rho >0,$ $k>0.$\bigskip\ 

The Lebesgue measure in $%
\mathbb{R}
^{d}$ has no atoms.Thus, the process is orderly.

\begin{theorem}
T\textit{he generalized Waring Process is a simple point process}
\end{theorem}

This follows directly from the Proposition 7.2.V, Daley \& Vere-Jones
(1988), since the generalized Waring Process in $%
\mathbb{R}
^{d}$ is orderly.

\subsection{Stationarity, $n$th-order Stationarity and Ergodicity}

The Lebesgue measure in $%
\mathbb{R}
^{d}$ is also is invariant under translations, hence the following can be
proved:

\begin{theorem}
\label{ths}Let $N\left( \cdot \right) $ be a generalized Waring process in $%
R^{d}$ with parameters $a>0$, $\rho >0,$ $k\in N$ .Then $N\left( \cdot
\right) $ is stationary.
\end{theorem}

\textit{Proof} We need to prove that for each $u\in R^{d}$ and all bounded
Borel sets $A\in \mathcal{B}_{%
\mathbb{R}
^{d}},$\ the avoidance function $P_{0}\left( \cdot \right) $ of the
generalized Waring process defined above satisfies $P_{0}\left( A\right) =$ $%
P_{0}\left( A+u\right) $ (see Daley \& Vere-Jones (1988), Theorem 10.1.III).

From the invariance of the Lebesgue measure on $%
\mathbb{R}
^{d}$ one can write $P_{0}\left( A\right) =$ $\dfrac{\rho _{\left( k\mu
\left( A\right) \right) }}{\left( \rho +a\right) _{\left( k\mu \left(
A\right) \right) }}=$ $\dfrac{\rho _{\left( k\mu \left( A+u\right) \right) }%
}{\left( \rho +a\right) _{\left( k\mu \left( A+u\right) \right) }}=$ $%
P_{0}\left( A+u\right) $ which proves the theorem.

A stationary Point process for which the $n$th-order moment measure exists
is $n$th-order stationary (see Daley \& Vere-Jones (1988)). Hence, using
Diggle \& Milne the following theorem and its corollary are trivial.

\begin{theorem}
The generalized Waring process in $R^{d}$ with parameters $a>0$, $\rho >n,$ $%
k\in N$ is $n$th-order stationary.
\end{theorem}

\begin{theorem}
The generalized Waring process in \bigskip $R^{d}$ is ergodic
\end{theorem}

\textit{Proof} \bigskip A necessary and sufficient criteria for a stationary
process to be ergodic is to be metrically transitive. From Lemma\ \ref{exl},
Appendix B, there exists one and only one root $x>0$ of the equation $\dfrac{%
\Gamma \left( \rho +x+a\right) }{\Gamma \left( \rho +x\right) }$ $=$\bigskip 
$b$ $>0.$ Let us consider $A$ a set in $R^{d}$ and let $S_{x\text{ }}$be the
shift operator. If $A$\ is such that $P\left( S_{x}A\cap A\right) $ $%
=P\left( A\right) $ then $\dfrac{\rho _{\left( k\mu \left( S_{x}A\cap
A\right) \right) }}{\left( \rho +a\right) _{\left( k\mu \left( S_{x}A\cap
A\right) \right) }}$ $=\dfrac{\rho _{\left( k\mu \left( A\right) \right) }}{%
\left( \rho +a\right) _{\left( k\mu \left( A\right) \right) }}.$ Hence we
obtain that $\ \dfrac{\Gamma \left( \rho +k\mu \left( S_{x}A\cap A\right)
+a\right) }{\Gamma \left( \rho +k\mu \left( S_{x}A\cap A\right) \right) }%
=\bigskip $ $\dfrac{\Gamma \left( \rho +k\mu \left( A\right) +a\right) }{%
\Gamma \left( \rho +k\mu \left( A\right) \right) }$and from Lemma\ \ref{exl}%
, Appendix B, follows that $\mu \left( S_{x}A\cap A\right) $ $=\mu \left(
A\right) .$ The last relation stands if $A=\Phi $ or $A=R^{d\text{ }}$which
does mean that $P\left( A\right) =0$ or $1.$This\ proves the theorem.

\section{Special Cases of the Generalized Waring Process}

In this section, we consider three instances of Genelarized Waring Processes
that may arise by multivariate extension, marginalization, projection, and
limiting arguments. Specifically, define the multivariate generalized Waring
process as a special case of the generalized Waring process on the product
space $S\times \left\{ 1,2,...m\right\} $ and show that marginals of a
multivariate GWP, as well as their sums, are all GWP as well. We then show
that GWP are closed under projection, and finally demonstrate how negative
binomial and Poisson processes can be seen as special (limiting) cases of
the GWP as some parameters are allowed to suitably diverge.

\subsection{The Multivariate Generalized Waring Process}

Consider the product space $S\times \left\{ 1,2,...m\right\} $ and let $%
\mathcal{B}_{S\times \left\{ 1,2,...m\right\} }$ be the associated product
Borel $\sigma $-algebra. Define the function $\nu :\mathcal{B}_{S\times
\left\{ 1,2,...m\right\} }\rightarrow R^{+\text{ }}$ such that for each $B=%
\overset{\infty }{\underset{i=1}{\sum }}A_{i}\times C_{i}\in \mathcal{B}%
_{S\times \left\{ 1,2,...m\right\} }(A_{i}\in \mathcal{B}_{S},$ $C_{i}\in 
\mathcal{P}\left( \left\{ 1,2,...m\right\} \right) ,$ $\nu \left( B\right) =%
\overset{\infty }{\underset{i=1}{\sum }}\mu \left( A_{i}\right) $ where $%
\mathcal{B}_{S}$ is the Borel $\sigma $-algebra and $\mu \left( \cdot
\right) $ some boundedly finite Borel measure. It is clear that $\nu \left(
\cdot \right) $ is a boundedly finite Borel measure on $S\times \left\{
1,2,...m\right\} .$ This allows us to define:

\textbf{Definition 2} The generalized Waring process with parameters $a$, $k$%
, $\rho $ and parameter measure $\nu \left( \cdot \right) $ on $S\times
\left\{ 1,2,...m\right\} $ is called the \textit{multivariate generalized
Waring process} with parameters $a$, $k$, $\rho $ and parameter measure $\mu
\left( \cdot \right) $ on $S$.

The multivariate GWP satisfies a number of convenient closure properties:

\begin{theorem}
Let ${\small N}\left( \cdot \right) $ be a multivariate generalized Waring
process with parameters $a$, $k$, $\rho $ and parameter measure $\mu \left(
\cdot \right) $ on $S.$ Then the following hold:

1. For every $i\in \left\{ 1,2,...m\right\} ,$ the marginal process ${\small %
N}_{i}\left( \cdot \right) =N\left( \cdot \times \left\{ i\right\} \right) $
is a GW process with parameters $a$, $k,$ $\rho $ and\textit{\ }parameter
measure\textit{\ }$\mu \left( \cdot \right) .$

2. $\overset{l}{\underset{j=1}{\sum }}N_{i_{j}}\left( \cdot \right) $ is a
generalized Waring process with parameters $a$, $\rho $ and\textit{\ }%
parameter measure\textit{\ }$kl\mu \left( \cdot \right) .$

3. For every finite collection of distinct indices $i_{1},i_{2},...,i_{l}$ $%
\in $ $\left\{ 1,2,...m\right\} $, $\{{\small N}_{i_{1}}\left( \cdot \right) 
{\small ,}$ ${\small N}_{i_{2}}\left( \cdot \right) {\small ,}$ ${\small ...,%
}$ ${\small N}_{i_{l}}\left( \cdot \right) \}$ is a multivariate generalized
Waring process with parameters $a,$ $\rho ,$ $k$ and parameter measure $\mu
\left( \cdot \right) .$

4. $\left\{ {\small N}_{i}\left( \cdot \right) {\small ,}\underset{j\neq i}{%
\sum }N_{j}\left( \cdot \right) \right\} $ is a bivariate generalized Waring
process with parameters $a$, $\rho $ and parameter measure $k\mu \left(
\cdot \right) ,\left( m-1\right) k\mu \left( \cdot \right) .$
\end{theorem}

\textit{Proof} For each bounded Borel set $A\in B_{s},$ the joint
distribution of $\{N_{1}\left( A\right) ,$ $N_{2}\left( A\right) ,$ $...,$ $%
N_{m}\left( A\right) \}$ is the $MGWD$ $(a;$ $k\mu \left( A_{1}\right) ,$ $%
k\mu \left( A_{2}\right) ,$ $...,$ $k\mu \left( A_{m}\right) ;$ $\rho )$.
From the structural properties of the multivariate generalized Waring
distribution (see Xekalaki (1986)), one has:

1. The distribution of $\left\{ {\small N}_{i}\left( A\right) {\small =x}%
_{i}\right\} $, for $i$ a given value on $\left\{ 1,2,...m\right\} $ is the
generalized Waring distribution with parameters $a$, $k\mu \left( A\right)
,\rho .$ By Theorem \ref{nsc}, this is a sufficient condition\ for the
process ${\small N}_{i}\left( \cdot \right) $ to be a generalized Waring
process.

2. The distribution of $\left\{ \overset{l}{\underset{j=1}{\sum }}%
N_{i_{j}}\left( A\right) {\small =x}_{i_{j}}\right\} ,$ is the generalized
Waring distribution with parameters $a$, $kl\mu \left( A\right) ,\rho .$ By
Theorem \ref{nsc}\ this is a sufficient condition\ for the process $\overset{%
l}{\underset{j=1}{\sum }}N_{i_{j}}\left( \cdot \right) $ to be a generalized
Waring process.

3. \bigskip\ For every $\left\{ A_{i_{1}},A_{i_{2}},...,A_{i_{l}}\in
B_{S}\right\} ,$ let us consider $\left\{ B_{1},B_{2},...,B_{m}\in
B_{S}\right\} $ where $B_{i}=B$ for $i\notin i_{1},i_{2},...,i_{l}$ and $%
B_{i}=A_{i_{j}}$ for $i=i_{j}.$ The joint distribution of $\{N_{1}\left(
B_{1}\right) ,$ $N_{2}\left( B_{2}\right) ,$ $...,$ $N_{m}\left(
B_{m}\right) \}$\ is the $MGWD$ $(a;$ $k\mu \left( B_{1}\right) ,$ $k\mu
\left( B_{2}\right) ,$ $...,$ $k\mu \left( B_{m}\right) ;$ $\rho ).$ From
the structural properties of the multivariate generalized Waring
distribution (see Xekalaki (1986)), it follows that the joint distribution
of $\{{\small N}_{i_{1}}\left( A_{i_{1}}\right) {\small ,}$ ${\small N}%
_{i_{2}}\left( A_{i_{2}}\right) {\small ,}$ ${\small ...,}$ ${\small N}%
_{i_{l}}\left( A_{i_{l}}\right) $ $\}$ is the $MGWD$ $(a;$ $k\mu \left(
A_{i_{1}}\right) ,$ $k\mu \left( A_{i_{2}}\right) ,$ $...,$ $k\mu \left(
A_{i_{l}}\right) ;$ $\rho )$ which proves part 3.

4. For every $A,B\in B_{S}$ let us consider $\left\{ B_{1},\text{ }B_{2},%
\text{ }...,B_{m}\text{ }\in B_{S}\right\} $ where $B_{i}=A$ and $B_{j}=B$
for $j\neq i.$ The joint distribution of $\{N_{1}\left( B_{1}\right) ,$ $%
N_{2}\left( B_{2}\right) ,$ $...,$ $N_{m}\left( B_{m}\right) \}$ is the 
\begin{equation*}
MGWD\left( a;\text{ }k_{1}\mu \left( B_{1}\right) ,\text{ }k_{2}\mu \left(
B_{2}\right) ,\text{ }...,\text{ }k_{m}\mu \left( B_{m}\right) ;\text{ }\rho
\right) .
\end{equation*}%
From the structural properties of the multivariate generalized Waring
distribution (see Xekalaki (1986)), it follows that the joint distribution
of $\{{\small N}_{i}\left( A\right) {\small ,}$ $\underset{j\neq i}{\sum }%
N_{j}\left( B\right) \}$ is the $BGWD$ $(a;$ $k\mu \left( A\right) ,$ $%
\left( m-1\right) k\mu \left( B\right) ;$ $\rho )$ which proves part
4.\bigskip

\subsection{Projections of Generalized Waring Processes}

Assume one has a product measurable space $(S_{1}\times S_{2},$ $\mathcal{B}%
_{S_{1}}\otimes \mathcal{B}_{S_{2}},$ $\mu _{1}\times \mu _{2})$ and let \ $%
N\left( \cdot \right) $ be a $GWP$ on that space, with parameteres $a,k,\rho 
$. Define $N_{S_{1}}\left( \cdot \right) $ and $N_{S_{2}}\left( \cdot
\right) $ to be the projections of $N\left( \cdot \right) $ onto $\left(
S_{1},\mathcal{B}_{S_{1}},\mu _{1}\right) $ and $\left( S_{2},\mathcal{B}%
_{S_{2}},\mu _{2}\right) $, respectively, defined by $N_{S_{1}}\left(
A\right) =N\left( A\times S_{1}\right) $ and $N_{S_{2}}\left( B\right)
=N\left( S_{2}\times B\right) .$ These projections will also be generalized
Waring processes:

\begin{theorem}
\label{projth} The projections $N_{S_{1}}\left( \cdot \right) $ an(d $%
N_{S_{2}}\left( \cdot \right) $ of a $GW\ $\ process $N\left( \cdot \right) $%
with parameteres $a,k,\rho ,$ $,$ onto the product measurable space $%
(S_{1}\times S_{2},$ $\mathcal{B}_{S_{1}}\otimes \mathcal{B}_{S_{2}},$ $\mu
_{1}\times \mu _{2})$ are also $GW$ processes with parameteres $a,b,\rho $
respectively onto $\left( S_{1},\mathcal{B}_{S_{1}},\mu _{1}\right) $ and $%
\left( S_{2},\mathcal{B}_{S_{2}},\mu _{2}\right) .$
\end{theorem}

\textit{Proof} \bigskip Let $\left\{ A_{i}\in \mathcal{B}%
_{S_{1}},i=1,2,...,l\right\} $ be finite family of disjoint bounded Borel
sets. The family $\{A_{i}\times S_{1}$ $\in \mathcal{B}_{S_{1}}\otimes 
\mathcal{B}_{S_{2}},$ $i=1,$ $2,$ $...,$ $l\}$ is also a finite family of
disjoint bounded sets. Hence, the joint distribution of $\left\{
N_{S_{1}}\left( A_{i}\right) =n_{i},i=1,...,l\right\} $ is the $MGWD$ $(a,$ $%
k\mu \left( A_{1}\right) ,$ $k\mu \left( A_{2}\right) ,$ $...,$ $k\mu \left(
A_{l}\right) ;$ $\rho )\ $which proves that the $N_{S_{1}}\left( \cdot
\right) $ is a $GWP$ with parameteres $a,k,\rho .$ The same argument yields
the result for $N_{S_{2}}\left( \cdot \right) .$

\subsection{The Poisson and the NB Processes as Limiting Cases of the
Generalized Waring Process}

We finally turn to demonstrate how negative binomial processes can be
obtained as limiting cases of the GWP. Doing so establishes that, even
though negative binomial processes cannot be orderly and stationary/ergodic
simultaneously, they can be approximated by a process with these properties.

\begin{theorem}
\label{gpPp} Let $N\left( \cdot \right) $ be a generalized Waring process
with parameters $a>0$, $\rho >0,$ and $k\in \emph{N}$, and with parameter
measure $\mu \left( \cdot \right) .$ Letting $k\rightarrow \infty $ and
setting $\rho =c\cdot k$ for $c>0$ a constant, the generalized Waring
process converges weakly to a Negative Binomial process with parameters $a$
and $c.$
\end{theorem}

\textit{Proof} Denote $N_{k}(\cdot )$, $k>0$ the generalized Waring process
indexed by the parameter $k$ and $N\left( \cdot \right) $ the Negative
Binomial process with parameters $a$ and $c.$ In order to prove that $%
N_{k}(\cdot )\underset{k\rightarrow \infty }{\rightarrow }N\left( \cdot
\right) $ weakly, it is sufficient to prove (see e.g. Daley \& Vere-Jones
(1988), Kallenberg (2002)):

(i). $P\left( N_{k}\left( A\right) =0\right) \underset{k\rightarrow \infty }{%
\rightarrow }P\left( N\left( A\right) =0\right) $ for all bounded $A$ of a
dissecting ring $\mathcal{T}$ of $S.$

(ii) That the generalized Waring process is uniformly tight.

In order to prove (i) we consider $P\left( N_{k}\left( A\right) =0\right) =%
\dfrac{\rho _{\left( k\mu \left( A\right) \right) }}{\left( \rho +\alpha
\right) _{\left( k\mu \left( A\right) \right) }}.$

We calculate:

$\dfrac{\rho _{\left( k\mu \left( A\right) \right) }}{\left( \rho +\alpha
\right) _{\left( k\mu \left( A\right) \right) }}=\dfrac{\rho _{\left(
a\right) }}{\left( \rho +k\mu \left( A\right) \right) _{\left( a\right) }}=%
\dfrac{ck_{\left( a\right) }}{\left( ck+k\mu \left( A\right) \right)
_{\left( a\right) }}$

$\cdot \dfrac{ck\left( ck+1\right) \cdot ...\cdot \left( ck+a-1\right) }{%
\left( ck+k\mu \left( A\right) \right) \left( ck+k\mu \left( A\right)
+1\right) \cdot ...\cdot \left( ck+k\mu \left( A\right) +a-1\right) }$

$=\dfrac{k^{a}c\left( c+\frac{1}{k}\right) \cdot ...\cdot \left( c+\frac{a-1%
}{k}\right) }{k^{a}\left( c+\mu \left( A\right) \right) \left( c+\mu \left(
A\right) +\frac{1}{k}\right) \cdot ...\cdot \left( c+\mu \left( A\right) +%
\frac{a-1}{k}\right) }$

$\underset{k\longrightarrow \infty }{\longrightarrow }\dfrac{c^{a}}{\left(
c+\mu \left( A\right) \right) ^{a}}=P\left( N\left( A\right) =0\right) $

To establish uniform tightness as required in (ii), we use two results
concerning regular and tight measures in a complete separable metric space $%
\mathcal{S}.$ A Borel measure is tight if and only if it is compact regular
(see e.g Lema A2.2.IV Daley \& Vere-Jones (1988)). In turn, a finite,
finitely additive, and nonnegative set function defined on the Borel sets of
a complete separable metric space $\mathcal{S}$ is compact regular if and
only if it is countably additive (see e.g Corollary A2.2.VII Daley \&
Vere-Jones (1988)). Therefore (ii) follows from the countable additivity
theorem (Theorem \ref{cath}), proven in earlier.

In turn, a Poisson process can be approximated by a negative binomial
process, so that it can also be approximated by a GWP:

\begin{theorem}
\label{pPpp}\bigskip Let $N\left( \cdot \right) $ be the limit proces $%
N\left( \cdot \right) $ $\ $of the previous Theorem, i.e. a Negative
Binomial process with parameters $a$ and $c.$ If $c\rightarrow \infty $ and $%
a=\lambda \cdot c$ where $\lambda >0$ is a constant, $N\left( \cdot \right) $
converges weakly to a Poisson process with parameter $\lambda $.
\end{theorem}

\textit{Proof} Write $\left\{ N_{c}(\cdot )\right\} ,$ $c>0$, to highlight
that the negative binomial process in question is indexed by $c$. We need to
show that there exists a Poisson process $M\left( \cdot \right) $ such that: 
$N_{c}(\cdot )\underset{c\rightarrow \infty }{\rightarrow }M\left( \cdot
\right) $ weakly. Following Daley \& Vere-Jones (1988), Lemma 9.I.IV, weak
convergence \ of the process and convergence of finite dimensional (fidi)
distributions are equivalent. So, in order to prove that $N_{c}(\cdot )%
\underset{c\rightarrow \infty }{\rightarrow }M\left( \cdot \right) $ weakly,
it is sufficient to prove that the fidi distributions of $N_{c}(\cdot )$
converge weakly to those of $M\left( \cdot \right) $ $.$

For every $\left\{ A_{_{1}},A_{2},...,A_{n}\in B_{S}\right\} $ we consider
the probability generating function $G_{c}$ $(A_{1},$ $A_{2},$ $...,$ $%
A_{n}; $ $z_{1},$ $...,$ $z_{n})$ of $N_{c}(\cdot )$ and obtain

$G_{n}$ $(A_{1},$ $A_{2},$ $...,$ $A_{n};$ $z_{1},$ $...,$ $z_{n})=\dfrac{1}{%
c}\left( c+\underset{i=1}{\overset{n}{\sum }}\left( 1-z_{i})\mu \left(
A_{i}\right) \right) \right) ^{-\lambda \cdot c}$

But $\dfrac{1}{c}\left( c+\underset{i=1}{\overset{n}{\sum }}\left(
1-z_{i})\mu \left( A_{i}\right) \right) \right) ^{-\lambda \cdot c}\underset{%
c\rightarrow \infty }{\rightarrow }\exp \left( -\lambda \underset{i=1}{%
\overset{n}{\sum }}\left( 1-zi)\mu \left( Ai\right) \right) \right) ,$ which
is the probability generating function $G$ $(A_{1},$ $A_{2},$ $...,$ $A_{n};$
$z_{1},$ $...,$ $z_{n})$ \ of the Poisson process with parameter $\lambda .$

\section{Discussion}

We have been able to define a new spatial point process for phenomena
characterised by over-dispersion, in great generality. The Generalized
Waring Process (GWP) has been shown to be able to simultaneously satisfy the
properties that negative binomial processes fail to (orderliness,
stationary, and ergodicity). Moreover, we have demonstrated that the new
process features appealing closure properties, in the sense that projection,
marginalization, and superposition all yield processes of the same GWD type,
with easily determinable parameters. By means of a a conditional property,
we have also illustrated that the process is straightforward to simulate.
These properties offer substantial advantages relative to existing
competitors of the negative binomial type, both from the theoretical and the
practical viewpoints, especially in terms of fitting the process on the
basis of a single realization. Indeed, we have shown that Generalized Waring
Point Process can even approximate negative binomial processes, giving a
positive resolution to the quandary posed in the conclusion of the paper by
Diggle and Milne: "Any view we adopt seems to seems to fall in a situation
from which progress looks difficult, and we conjecture that no stationary,
ergodic, orderly negative binomial processes exist." Though negative
binomial processes may fail to simultaneously verify
orderliness/stationarity/ergodicity, they can be well approximated by
flexible and tractable processes of the GWP class that do verify these
properties.

Potential further advantages of the generalized Waring Process relative to
negative binomial\ processes\ may arise in the context of compounding (or
clustering) and mixing (or heterogeneity). In particular, Cane (1974,1977)
has demonstrated that one cannot distinguish between compounding and
heterogeneity under a negative binomial distribution: given a total of $n$
events, the distribution of event times is the same, whether one the model
arose out of mixing or compounding. In contrast, Xekalaki (1983b)
demonstrated that discriminating between clustering and mixing may well be
possible in the context of the Generalized Waring Distribution, by showing
that the conditional distribution of the times of events given their total
is different under compounding and under mixing (see also Xekalaki (2006,
2014, 2015). This property can then be used in order to distinguish
clustering, which may otherwise be confounded with compounding.

\section{References}

Ajiferuke, I.,Wolfram, D. \& Xie, H. (2004). Modelling website visitation
and resource usage characteristics by IP address data. In: \textit{%
Proceedings of the 32nd Annual Conference of the Canadian Association for
Information Science (2004) (eds H. Julien and S. Thompson)}, Available at:
www.cais-acsi.ca/proceedings/ 2004/ajiferuke\_2004.pdf

Bai-ni Guo, Ying-jie Zhang \& Feng Qi. (2008). Refinements and Sharpenings
of Some double Inequalities for Bounding the Gamma Function,\textit{\ J.
Inequal. Pure Apl. Math.}, \textbf{9} , Issue 1, Article 17.

Bates, G.E. (1955). Joint distributions of time intervals for the ovurrence
of successive accidents in a generalized Polya scheme. \textit{Annals of
Math. Statist.} \textbf{26}, 705-720.

Boswell, M.T. \& Patil, G.P. (1977). Chance mechanisms generating the
negative binomial distribution, in: \textit{Random Counts in Scientific Work
1}: \textit{Random Counts in Models and Structures} (The Pennsylvania State
Univ. Press) 3-22.

Burnett, R.T. \& Wasan, M.T. (1980). The negative binomial point process and
its inference, in: \textit{Multivar. Statist. Analysis} (eds R.P. Gupta),
31-45, North-Holland Publishing Co., Amsterdam-New York.

Cane, V. R. (1974). The concepts of accident proneness. \textit{Bull. of
Institute of Math., Bulgarian Academy of Sci.} \textbf{15}, 183-189.

Cane, V.R. (1977). A class of non-identifiable stochastic models. \textit{J.
Appl. Probab}. \textbf{14}, 475-782

Chatfield, C. \& Theobald, C.M. (1973). Mixtures and Random Sums, \textit{%
The Statistician}, \textbf{22(3)}, 281-287.

Chetwynd, A. G. \& Diggle, P. J. (1998). On estimating the Reduced Second
Moment Measure of a Stationary spatial Point Process, \textit{Austral. \&
New Zealand J.Statist}. \textbf{40(1)}, 11-15.

Cliff A.D. \& Ord J.K. (1973). Spatial Autocorrelation (Pion, London).

Cox, D.R. \& Isham, V.I. (1980). Point Processes. Chapman \& Hall. ISBN
0-412-21910-7.

Cox, D.R. \& Lewis, P.A.W. (1966). The Statistical Analysis of Series of
Events. London: Methuen.

Cresswell, W. L. \& Froggatt, P. (1963). The Causation of Bus Driver
Accidents. London. Oxford University Press.

Daley, D.J. \& Vere-Jones, D. (1972). A summary of the mathematical theory
of point processes. In: \textit{Stoch. Point Proc.} (eds P.A.W. Lewis).
299-383.

Daley, D.J. \& Vere-Jones, D. (1988). An Introduction to the Theory of Point
Processes, Springer

Diggle, P.J. \& Chetwynd, A. G. (1991). Second- order Analysis of Spatial
Clustering for Inhomogeneous Populations. \textit{Biometrics}, \textbf{47},
1155-1163

Diggle, P.J. \& Milne R.K. (1983). Negative Binomial Quadrat Counts and
Point Processes, \textit{Scandinavian J. of Statist.} Vol. \textbf{10(4)},
257-267.

Feller, W. (1968). An introduction to probability theory and its
applications. Vol. I, 3rd ed. Wiley, New York.

Fisher, L. (1972). A survey of the mathematical theory of multidimensional
point processes. In: \textit{Stoch. Point Proc.: Stat. Anal., Theory, and
Applications} (eds P. A. W. Lewis), 468-51, New York: Wiley-Interscience.

Gregoire, G. (1983). Negative binomial distributions for point processes, 
\textit{Stoch. Proc. Appl}. \textbf{16}, 179-188.

Grandell, J. (1997). Mixed Poisson Processes. Chapman \& Hall, London

Gurland, J. (1959). Some applications of the negative binomial and other
contagious distributions. \textit{American J. of Public Health} \textbf{49},
1388-1399.

Hjort, N. (1990). Nonparametric Bayes estimators based on beta processes in
models for life history data. \textit{Annals of Statist.} \textbf{18},
1259--1294.

Irwin, J. O. (1941). Discussion on Chambers and Yule's paper. \textit{J.
Roy. Statist. Soc}.\textit{, }Supplement\textbf{\ 7}, 101-109.

Irwin, J.O. (1963). The place of mathematics in medical and biological
statistics. \textit{J. Roy. Statist. Soc}., (\textbf{A)126}, 1-44.

Irwin, J.O. (1968).The Generalized Waring Distribution. Applied to Accident
Theory, \textit{J. Roy. Statist. Soc}., \textbf{(A)131}, 205-225.

Irwin, J.O. (1975). The Generalized Waring Distribution, \textit{JJ. Roy.
Statist. Soc.}, \textbf{(A)138}, 18-31 (Part I), 204-227 (Part II), 374-384
(Part III)

Janardan, K. G. \& Patil, G. P. (1972). A unified approach for a class of
multivariate hypergeometric models. \textit{Sankhya - Ind. J. Stat. } 
\textit{(A)34}, 363-376.

Johnson, N. L., Kotz, S. \& Kemp, A. W. (1992). Univariate Discrete
Distributions. John Wiley and Sons, Inc., New York.

Kallenberg, O. (2002). Foundations of Modern Probability, Springer.

Kemp, C. D. (1970). 'Accident proneness' and discrete distribution theory.
In: \textit{Random counts in Scientific Work}, (eds G.P. Patil) \textbf{2},
41-65. State College: Pennsilvania State University Press.

Kemp, A. W. \& Kemp, C. D. (1975). Models for Gaussian hypergeometric
distribution. In: \textit{Statistical Distributions in Scientific Work.
Dodrecht: Reidel} (eds G. P. Patil, S. Kotz \& J. K. Ord) \textbf{1}, 31-40.

Matern, B. (1971). Doubly stochastic Poisson processes in the plane. In: 
\textit{Statistical Ecology}, (eds. G.P. Patil, E.C. Pielou \& W.E. Waters),
195--213. University Park: Pennsylvania.

Merkle M. (2001). Conditions for Convexity of a Derivative and Applications
to the Gamma and Digamma Function, \textit{Serb. Math. Inform.} \textbf{16},
13-20.

Muller, J. \& Waagepetersen, R.P. (2004). Statistical Inference and
Simulation for Spatial Point Processes. Chapman and Hall/CRC, Boca Raton.

Newbold, E. M. (1927). Practical Applications of the Statistics of Repeated
Events, Particularly to Industrial Accidents. \textit{J. Roy. Statist. Soc}%
., \textbf{90}, 487-547.

Ripley, B. D. (1988). Statistical Inference for Spatial Processes, Cambridge
University Press

Ross, S.M. (1995). Stochastic Processes. Wiley. ISBN 978-0471120629

Ross, S.M. (2007). Introduction to probability models Amsterdam; Boston :
Elsevier/Academic Press, 9th ed.

Ramakrishnan, A. (1951). Some simple stochastic processes, \textit{J. Roy.
Statist. Soc}., \textbf{(B)13(1)}, 131-140.

Shaw, L. \& Sichel, H. S. (1971). Accident Proneness. Oxford: Pergamon Press.

Sibuya, M. (1979). Generalized hypergeometric, digamma and trigamma
distributions. \textit{Annals Inst. Statistical Math.}, \textbf{31/A},
373--390.

Sibuya, M. \& Shimizu, R. (1981). Classification of the generalized
hypergeometric family of distributions. \textit{KEIO Sci. and Techn. Reports}%
, \textbf{34}, 1--38.

Snyder, D.L. \& Miller, M.I. (1991) Random Point Processes in Time and
Space. Springer-Verlag. ISBN 0-387-97577-2.

Xekalaki, E. (1981). Chance Mechanisms for the Univariate Generalized Waring
Distribution and Related Characterizations. \textit{Statist. Distr. in Sci.
Work}, \textbf{4}, 157-172.

Xekalaki, E.(1983a). Infinite Divisibility, Completeness and Regression
Properties of the Univariate Generalized Waring Distribution, \textit{Annals
Inst. Statist. Math.,} \textbf{(A) 35}, 161-171.

Xekalaki, E. (1983b). The Univariate Generalized Waring Distribution in
Relation to Accident Theory. Proneness, Spells or contagion?, \textit{%
Biometrics} \textbf{39 (3)}, 39-47.

Xekalaki, E. (1983c). A Property of the Yule Distribution and its
Applications. \textit{Commun.Statist.-Theor.Meth.}, \textbf{(A)12(10)},
1181-1189.

Xekalaki, E.(1983d) Hazard Functions and Life Distributions in Discrete
Time. \textit{Commun.Statist.-Theor.Meth.}, \textbf{12 (21)}, 2503-2509.

Xekalaki E. \&\ Panaretos J. (1983). Identifiability of Compound Poisson
Distributions\textquotedblright . \textit{Scand. Actuarial J.}, \textbf{66},
39-45.

Xekalaki, E. (1984). The Bivariate Generalized Waring Distribution and its
Application to Accident Theory, \textit{J. Roy. Statist. Soc}., \textbf{%
(A)47(3)}, 448-498.

Xekalaki, E. (1985a). Factorial Moment Estimation for the Bivariate
Generalized Waring Distribution. \textit{Statistiche Hefte }\textbf{26},
115-129.

Xekalaki, E.(1985b). Some Identifiability Problems Involving Generalized
Waring Distributions. \textit{Publicationes Mathematicae}, \textbf{32},
1985, 75-84.

Xekalaki, E. (1986). The Multivariate Generalized Waring Distribution. 
\textit{Commun.Statist.-Theor.Meth.},\textbf{15(3)} 1047-1064.

Xekalaki, E. (2006). Under- and Overdispersion. In: \textit{Encyclopedia of
Actuarial Science}, 3, 1700-1705. Available at DOI:
10.1002/9780470012505.tau003

Xekalaki, E. \& Zografi, M. (2008). The Generalized Waring Process and its
Applications -- Commun.Statist.-Theor.Meth., \textbf{37(12)}, 1835-1854

Xekalaki, E. (2014) On The Distribution Theory of Over-Dispersion. \textit{%
J. Statist. Distr. \& Applications}, \textbf{1(19)}, on line version DOI:
10.1186/s40488-014-0019-z

Xekalaki, E. (2015). Under- and Overdispersion. Wiley StatsRef: Statistics
Reference Online. 1--9. in: \textit{Encyclopedia of Actuarial Science}.
(Update based on the 2006 original article by E. Xekalaki) Available at:
http://onlinelibrary.wiley. com/doi/10.1002/9781118445112.stat04407.pub2.

Westcott, M. (1971). On existence and mixing results for cluster point
processes. \textit{J. Roy. Statist. Soc.} \textbf{(B)33}, 290-300.

Westcott, M. (1972). The probability generating functional. \textit{J.
Austral. Math. Soc.} \textbf{14}, 448-466.

Willmot, G. E. (1986) Mixed compound distributions. \textit{ASTIN Bulletin} 
\textbf{16/S}, 59-S79.

Zografi, M. \& Xekalaki, E. (2001). The Generalized Waring Process. In: 
\textit{Proceedings of the 5th Hellenic-European Conference on Computer
Mathematics and its Applications}, (eds E.A. Lypitakis), 886-893, Athens,
Greece.\bigskip

\bigskip

\textbf{Corresponding author's address}: Professor Evdokia Xekalaki,
Department of Statistics, Athens University of Economics, 76 Patision St.,
10434 Athens, Greece

\bigskip

\bigskip

\bigskip

\bigskip

\bigskip

\section{\protect\bigskip Appendix}

\subsection{Existence Lemma}

\bigskip The following Lemma proves that the equation $P_{0}\left( A\right) =%
\dfrac{\rho _{\left( k\mu \left( A\right) \right) }}{\left( \rho +a\right)
_{\left( k\mu \left( A\right) \right) }}$ always has a unique solution. This
result is used in the proof of existence of GWP.

\begin{lemma}
\label{exl}\ For each $a>0$, $\rho >0,0\leq P_{0}\leq 1,$ there exists one
and only one root $x>0$ of the equation $\frac{\Gamma \left( \rho
+x+a\right) }{\Gamma \left( \rho +x\right) }=$\bigskip $\frac{\Gamma \left(
\rho +a\right) }{P_{0}\Gamma \left( \rho \right) }$
\end{lemma}

\bigskip Proof

It has been proved (see Bai-ni, Ying-jie and Feng \cite{2} Theorem 3\textbf{%
) }that $\dfrac{\Gamma \left( y\right) }{\Gamma \left( x\right) }>\dfrac{%
y^{y-\gamma }}{x^{x-\gamma }}e^{x-y}$ for all for $y>x\geq 1$, where\ $%
\gamma $ stands for the Euler-Mascheroni constant

\bigskip Hence we can obtain%
\begin{equation*}
\frac{\Gamma \left( \rho +x+a\right) }{\Gamma \left( \rho +x\right) }>\frac{1%
}{e^{a}}\frac{\left( \rho +x+a\right) ^{\rho +x+a-\gamma }}{\left( \rho
+x\right) ^{\rho +x-\gamma }}
\end{equation*}

\bigskip Now consider the function%
\begin{equation*}
f\left( x\right) =\frac{\left( \rho +x+a\right) ^{\rho +x+a-\gamma }}{\left(
\rho +x\right) ^{\rho +x-\gamma }}e^{-a}
\end{equation*}

\bigskip with derivative%
\begin{eqnarray*}
f^{\prime }\left( x\right) &=&e^{-a}\frac{\left[ \ln \left( \rho +x+a\right)
+1-\frac{\gamma }{\rho +x+a}\right] \left( \rho +x+a\right) ^{\rho
+x+a-\gamma }}{\left( \rho +x\right) ^{2\left( \rho +x-\gamma \right) }} \\
&&e^{-a}\frac{-\left[ \ln \left( \rho +x\right) +1-\frac{\gamma }{\rho +x}%
\right] \left( \rho +x\right) ^{\rho +x-\gamma }}{\left( \rho +x\right)
^{2\left( \rho +x-\gamma \right) }}
\end{eqnarray*}

\bigskip \bigskip The functions $\varphi \left( x\right) =\left( 1+\ln x-%
\frac{\gamma }{x}\right) $ and $\omega \left( x\right) =x^{x}$ are
increasing for $x>1$, since $\varphi ^{\prime }\left( x\right) =\left( \frac{%
1}{x}+\frac{\gamma }{x^{2}}\right) >0$ if $x>-\gamma ,$ \bigskip \bigskip $%
\omega ^{\prime }\left( x\right) =\left( 1+\ln x\right) x^{x}>0$ if $x>1.$
Hence the function $g\left( x\right) =\varphi \left( x\right) \omega \left(
x\right) =\left( 1+\ln x-\dfrac{\gamma }{x}\right) x^{x}$ is increasing for $%
x>1.$

Therefore, $\left[ \ln \left( \rho +x+a\right) +1-\dfrac{\gamma }{\rho +x+a}%
\right] $ $\left( \rho +x+a\right) ^{\rho +x+a-\gamma }$ $\ -$ $\ [\ln
\left( \rho +x\right) +1-\dfrac{\gamma }{\rho +x}]$ $\left( \rho +x\right)
^{\rho +x-\gamma }$ $\ >0$ for $x>1$ so that $f^{\prime }\left( x\right) >0$%
, proving that $f\left( x\right) $ is increasing for $x>1$.

\bigskip In summary, we can state that $\forall b\in R,$ $\exists x>1,$ such
that $f\left( x\right) >b.$

$\bigskip $Let us now consider $b=\dfrac{\Gamma \left( \rho +a\right) }{%
P_{0}\Gamma \left( \rho \right) }$. For that value $\exists x>0,$ such that $%
f\left( x\right) >b$. Clearly for $x=0,$ $\dfrac{\Gamma \left( \rho
+x+a\right) }{\Gamma \left( \rho +x\right) }=\dfrac{\Gamma \left( \rho
+a\right) }{\Gamma \left( \rho \right) }<b.$ The function $\dfrac{\Gamma
\left( \rho +x+a\right) }{\Gamma \left( \rho +x\right) }$ is continuous for $%
x>0$ as a ratio of two continuous functions. So, applying Bolzano's Theorem
to $\dfrac{\Gamma \left( \rho +x+a\right) }{\Gamma \left( \rho +x\right) }%
-b, $ we obtain that $\exists x>0$ such that $\dfrac{\Gamma \left( \rho
+x+a\right) }{\Gamma \left( \rho +x\right) }-b=0.$

$\bigskip $On the other hand%
\begin{equation*}
\dfrac{d}{dx}\left[ \dfrac{\Gamma \left( \rho +x+a\right) }{\Gamma \left(
\rho +x\right) }\right] =\dfrac{d}{dx}\left[ \exp \left( \ln \dfrac{\Gamma
\left( \rho +x+a\right) }{\Gamma \left( \rho +x\right) }\right) \right]
\end{equation*}

\begin{equation*}
=\dfrac{\Gamma \left( \rho +x+a\right) }{\Gamma \left( \rho +x\right) }%
\dfrac{d}{dx}\left[ \ln \Gamma \left( \rho +x+a\right) -\ln \Gamma \left(
\rho +x\right) \right]
\end{equation*}

\begin{equation*}
=\dfrac{\Gamma \left( \rho +x+a\right) }{\Gamma \left( \rho +x\right) }\left[
\Psi \left( \rho +x+a\right) -\Psi \left( \rho +x\right) \right]
\end{equation*}%
and using the relation $\Psi \left( t\right) =-\gamma +\underset{i=0}{%
\overset{\infty }{\sum }}\left( \dfrac{1}{i+1}-\dfrac{1}{i+t}\right) $ where 
$\gamma $\ is the Euler-Mascheroni constant, we obtain%
\begin{eqnarray*}
\dfrac{d}{dx}\left[ \dfrac{\Gamma \left( \rho +x+a\right) }{\Gamma \left(
\rho +x\right) }\right] &=&\dfrac{\Gamma \left( \rho +x+a\right) }{\Gamma
\left( \rho +x\right) } \\
\overset{\infty }{\sum_{i=0}}\left( \frac{1}{i+\rho +x}-\frac{1}{i+\rho +x+a}%
\right) &>&0
\end{eqnarray*}%
which proves the Lemma.

\end{document}